\documentclass[journal]{IEEEtran}
\usepackage{multirow}
\usepackage{multicol}
\usepackage{lipsum}
\usepackage{graphicx}
\usepackage{graphics}
\usepackage{amsmath}
\usepackage{amsbsy}
\usepackage{blindtext}
\usepackage{tcolorbox}
\usepackage{amssymb}
\usepackage{scalerel}
\usepackage{mathabx}
\usepackage{verbatim}
\usepackage{booktabs}
\usepackage{rotating,tabularx}
\usepackage{mathrsfs} 
\usepackage{url}
\usepackage{adjustbox}
\usepackage{soul}
\usepackage{xcolor}
\usepackage{cite}
\usepackage{tikz}
\usepackage{color, colortbl}
\usepackage[first=0,last=9]{lcg}
\definecolor{LightGray}{rgb}{0.7,0.7,0.7}

\usepackage[wby]{callouts}
\usetikzlibrary{shapes.multipart}
\usepackage[compact]{titlesec}         
\AtBeginDocument{
  \setlength\abovedisplayskip{0pt}
  \setlength\belowdisplayskip{0pt}}
\usepackage{array}
\usepackage{makecell}

\allowdisplaybreaks

\usepackage{amsthm}
\theoremstyle{definition}

\theoremstyle{remark}

\usepackage[utf8]{inputenc}
\usepackage[english]{babel}
\usepackage{float}

\usepackage[caption=false,font=footnotesize]{subfig}
\usepackage[T1]{fontenc}
\usepackage{wrapfig}
\usepackage{scalerel,stackengine}
\newcommand\reallywidecheck[1]{%
\savestack{\tmpbox}{\stretchto{%
  \scaleto{%
    \scalerel*[\widthof{\ensuremath{#1}}]{\kern-.6pt\bigwedge\kern-.6pt}%
    {\rule[-\textheight/2]{1ex}{\textheight}}
  }{\textheight}%
}{0.5ex}}%
\stackon[1pt]{#1}{\scalebox{-1}{\tmpbox}}%
}
\stackMath
\IEEEoverridecommandlockouts
\IEEEoverridecommandlockouts
\renewcommand\Re{\operatorname{Re}}
\renewcommand\Im{\operatorname{Im}}

\newif\ifarxiv
 \arxivtrue

\title{Assessing the impact of Higher Order Network Structure on Tightness of OPF Relaxation }

\ifarxiv
\author{Nafis Sadik,$^{\ast}$ and Mohammad Rasoul Narimani,$^{\dagger}$%
\thanks{This work was supported by NSF under Award Number 2308498.}
\thanks{${\ast}$: College of Engineering and Computer Science, Arkansas State University. nafis.sadik@smail.astate.edu.}%
\thanks{$^{\dagger}$: Department of Electrical and Computer Engineering, California State University Northridge (CSUN). Rasoul.narimani@csun.edu.}}%

\else
\author{Mohammad Rasoul Narimani and Nafis Sadik}%
\fi

\begin{document}
\maketitle

\begin{abstract}
AC optimal power flow (AC OPF) is a fundamental problem in power system operation and control. Accurately modeling the network physics via the AC power flow equations makes AC OPF a challenging nonconvex problem that results in significant computational challenges. To search for global optima, recent research has developed a variety of convex relaxations to bound the optimal objective values of AC OPF problems. However, the quality of these bounds varies for different test cases, suggesting that OPF problems exhibit a range of difficulties. Understanding this range of difficulty is helpful for improving relaxation algorithms.
Power grids are naturally represented as graphs, with buses as nodes and power lines as edges. Graph theory offers various methods to measure power grid graphs, enabling researchers to characterize system structure and optimize algorithms. Leveraging graph theory-based algorithms, this paper presents an empirical study aiming to find correlations between optimality gaps and local structures in the underlying test case's graph. Network graphlets, which are induced subgraphs of a network, are used to investigate the correlation between power system topology and OPF relaxation tightness. Specifically, this paper examines how the existence of particular graphlets that are either too frequent or infrequent in the power system graph affects the tightness of the OPF convex relaxation.
Numerous test cases are analyzed from a local structural perspective to establish a correlation between their topology and their OPF convex relaxation tightness.


\end{abstract}


\IEEEpeerreviewmaketitle

\section{Introduction}

\IEEEPARstart{T}{he} optimal power flow (OPF) problem seeks an operating point that optimizes a specified objective function (often generation cost minimization) subject to constraints from the network physics and engineering limits. Using the nonlinear AC power flow model to accurately represent the power flow physics results in the AC OPF problem, which is non-convex,
generally NP-Hard problems~\cite{bienstock2015nphard} and may have local optima~\cite{low2014convex}. The inclusion of AC power flow in the OPF problem presents non-convex feasible spaces, leading to considerable computational complexity~\cite{narimani2018empirical, narimani2020strengthening}. Many convex relaxation techniques have been used to solve OPF problems~\cite{convcoff}. These relaxation techniques converge to a lower bound for the OPF problem. In some cases if the relaxation solution satisfies specific condition, the calculated lower bound by relaxation method can be inferred as the global solution of the problem. An optimality gap can be referred as the percent difference between the objective values for a local solution and the lower bound obtained by  relaxation techniques. A relatively smaller optimality gap certifies global optimal solution of corresponding local solution~\cite{barzegar2019method}. 

Many research efforts have been done in the past decades to develop convex relaxation algorithms for NP-hard, non-convex problems~\cite{gopinath2020proving}. Numerous relaxation techniques including Semi-definite programming (SDP)~\cite{Lavaei}, second order cone programming (SOCP)~\cite{Bose}, quadratic convex relaxation (QC)~\cite{coffrin2015qc}. QC relaxation encloses the trigonometric, squared, and product terms in a polar representation of power flow equations within convex envelopes~\cite{coffrin2015qc}. We utilized the QC relaxation in this paper for investigating correlation between relaxation optimality gap and topology of the test cases.  

The correlation between power systems and graph theory is strong, as graph theory provides a powerful tool for analyzing and optimizing complex power systems. This correlation was first explored in~\cite{seshu1961linear}. Graph theoretical analysis approach have been used in many power system applications such as system vulnerability~\cite{wang2010node,cuadra2015critical}, detecting structural anomalies~\cite{cuffe2015assortativity}, identifying critical components in power systems~\cite{Narimani1,Narimani2, Boyaci_propagation}, and generating authentic synthetic grids~\cite{birchfield2016grid}. In the mean time, some research has been done on optimal power flow problem with the perspective from network science. According to~\cite{chitra2015optimal}, power flow in power networks can be traced by using graph theoretical algorithms such as breadth first search and depth first search. Connection between cliques and semi-definite solver performance is discussed in~\cite{kersulis2018topological} and it is suggested that semi-definite constraints can be decomposed into smaller constraints according to maximal cliques of the power network. 
In addition, to make the SDP relaxation problem easier to solve, identifying the problematic lines that contribute to its computational complexity is crucial. This can be accomplished by using graph theory techniques, such as tree decomposition, to analyze the structure of the graph and identify those lines, as described in~\cite{madani2015promises}.
Though not focused in optimal power flow, the relationship between the network topology, as characterized by the maximal cliques, and the number of power flow solution has been explored in~\cite{molzahn2016toward}. More recently, graph neural networks are gaining attention to solve OPF problems. Particularly, graph neural networks can be used to approximate optimal interior point optimizer solution in OPF~\cite{owerko2020optimal}. However, to the best of our knowledge, no studies have investigated the correlation between the optimality gap of OPF relaxations and the local structures of the underlying power grid's graph in power system test cases.
%

Graphlets and motifs are essential tools in complex network analysis. They are subgraphs that occur frequently in a given network and can provide insights into the network's structure and function. Graphlets are small, connected subgraphs that can be used to characterize the local structure of a network. Motifs are larger subgraphs that occur more frequently than expected by chance and can represent functional units in the network. By identifying and analyzing graphlets and motifs, researchers can gain a deeper understanding of the network's organization and dynamics, and identify key nodes or pathways that are critical to the network's function. Graphlet and motif analysis can be used to identify similarities and differences, and classify networks based on their structure and function. 
Use of network motifs statistically to extract information about the local structure of the data was first proposed in~\cite{milo2002network}. At first, analysis of motif detection in a network was exclusively restricted in the field of bio-informatics~\cite{lacroix2005reaction}~\cite{prill2005dynamic}. Using network motifs, comparative grid vulnerability analysis in case of contingency has been discussed in~\cite{dey2017motif} and it was found that vulnerable and robust grid have different motif patterns and decay of motifs also reveals a different pattern when comparing robust and fragile grids. What these local structure of a grid network means to a power grid in case of a contingency has already been discussed~\cite{abedijaberi2018motif} and it has been found that certain motifs did play determining the robustness of a network. 

This paper first explores network graphlet patterns across various test cases and next determines the optimality gaps for QC relaxation in those OPF test cases. In particular, we attempt to ascertain if there is any connection between the local structure of the network and the subsequent OPF optimality gap for that specific network. To be more specific, this paper delves deeper to discover if any graphlets contribute to a larger optimality gap in OPF test cases. By identifying graphlets that contribute to a larger optimality gap in test cases, we can identify important nodes in test systems, wherein enforcing redundant constraints on them reduces the optimality gap for those test cases. Moreover, such analysis can lead to the development of more effective optimization algorithms tailored to specific network topologies and requirements. Thus, identifying significant graphlets that contribute to larger optimality gaps can ultimately lead to the development of more efficient and reliable network infrastructures.

This paper is organized as follows. Sections~\ref{Overview of the Optimal Power Flow Problem} and~\ref{Convex Relaxations in Optimal Power flow} review the OPF formulation and the QC relaxation, respectively. Section~\ref{Network Motifs} describes the network graphlet theory. Section~\ref{Statistical Analysis} then presents an algorithm by which graphlets will be detected. Section~\ref{Results} discusses the numerical result and Section~\ref{Conclusion} concludes the paper.

\section{Overview of the Optimal Power Flow Problem}
\label{Overview of the Optimal Power Flow Problem}
This section formulates the OPF problem using a polar
voltage phasor representation. 
Let's assume voltage at bus $i$ and $j$ are $V_i$ and $V_j$, respectively. Current flowing from bus $i$ to bus $j$ is $I_{ij}$, AC power generation of bus $i$ is $S_i^g$, complex power demand of bus $i$ is $S_i^d$ and complex power flowing from bus $i$ to bus $j$ is $S_{ij}$. $E$ and $F$ are the set of sending and receiving ends of lines, i.e. edges, respectively. $Y_{i j}$ is the admittance of line $i$ to $j$.  Power flow equation for all buses can be written as follows.

\begin{subequations}
\label{OPF formulation}
\begin{align}
\label{eq:objective}
S_i^g-S_i^d = \sum_{\substack{(i,j)\in \ {E \cup F}}} \!S_{ij}\quad \forall i \in \mathcal{\text{N}}\\
S_{i j}=Y_{i j}^{*} V_{i} V_{i}^{*}-Y_{i j}^{*} V_{i} V_{j}^{*}\quad (i, j) \in E \cup F
\end{align}
\end{subequations}

The OPF problem consists of different engineering constraints that should be enforced along with the power flow equations. Generators in the system should produce active and reactive power within their limits which can be addressed by following constraints. 

\begin{align}
\label{eq:objective2}
S_i^{gl} \le S_i^g \le S_i^{gu} \quad \forall i \in \mathcal{\text{N}}
\end{align}

Line thermal limit is another constraint that enforces an upper bound on apparent power flow in lines. 

\begin{align}
\label{eq:objective2}
|S_{ij}| \le s_{ij}^u  \quad\forall i \in \mathcal{\text{N}}
\end{align}

Bus voltage limits of any grid are defined by national grid code of any country and it is typically $\pm10\%$ of nominal grid voltage~\cite{kundur2004definition}.

\begin{align}
\label{eq:objective2}
v_i^{l} \le |V_i| \le v_i^{u} \quad\forall i \in \mathcal{\text{N}}
\end{align}

For ease of power flow formulation squaring this equation gives us,

\begin{align}
\label{eq:objective2}
{v_i^{l}}^2 \le {|V_i|}^2 \le {v_i^{u}}^2\quad\forall i \in \mathcal{\text{N}}
\end{align}

For power flow between buses, voltage angle difference between buses must be confined. Accordingly, voltage angle difference upper and lower limits can be shown as,

\begin{align}
\label{eq:angle difference}
-\theta_{ij}^{\Delta} \le \angle(V_{i} V_{j}^{*}) \le \theta_{ij}^{\Delta}\quad\forall {(i,j)} \in \mathcal{\text{E}}
\end{align}

To convexify the OPF problem, phase angle difference should be limited within $[0,\pi/2]$~\cite{coffrin2014linear}. 

\begin{align}
\label{eq:theta pi/2}
0 \le \theta_{ij}^{\Delta} \le {\pi}/2\quad\forall {(i,j)} \in \mathcal{\text{N}}
\end{align}

If we observe equation \eqref{eq:angle difference} as linear relation of real and imaginary parts of $\angle(V_{i} V_{j}^{*})$ then it can be shown that,

\label{eq:objective2}
\begin{align}
\label{eq:objective2}
\tan \left(-\theta_{i j}^{\Delta}\right) \Re\left(V_{i} V_{j}^{*}\right) \leqslant \Im\left(V_{i} V_{j}^{*}\right) \leqslant \tan \left(\theta_{i j}^{\Delta}\right) \Re\left(V_{i} V_{j}^{*}\right)
\end{align}

Objective of the OPF problem is to minimize generator fuel costs and that can be defined as,

\label{eq:objective2}
\begin{align}
\label{eq:objective2}
\min \sum_{i \in \mathcal{N}} c_{2, i}\left(P_{i}^{g}\right)^{2}+c_{1, i} P_{i}^{g}+c_{0, i}
\end{align}

In overall, the OPF problem can be written as,

\begin{subequations}
\label{OPF Formulation overall}
\begin{align}
\label{eq:objective}
\text{variables}: S_{i}^{g}(\forall i \in N), V_{i}\quad(\forall i \in N)
\end{align}
\begin{align}
\min \sum_{i \in \mathcal{N}} c_{2, i}\left(P_{i}^{g}\right)^{2}+c_{1, i} P_{i}^{g}+c_{0, i}
\end{align}
\begin{align}
\text{subject to}:v_{i}^{l}
\leqslant\left|V_{i}\right| \leqslant v_{i}^{u}\quad \forall i \in N
\end{align}
\begin{align}
S_{i}^{g l} \leqslant S_{i}^{g} \leqslant S_{i}^{g u}\quad \forall i \in N 
\end{align}
\begin{align}
\left|S_{i j}\right| \leqslant s_{i j}^{u}\quad \forall(i, j) \in E \cup F 
\end{align}
\begin{align}
S_{i}^{g}-S_{i}^{d}=\sum_{(i, j) \in E \cup F} S_{i j}\quad \forall i \in N 
\end{align}
\begin{align}
\label{10g}
S_{i j}=Y_{i j}^{*} V_{i} V_{i}^{*}-Y_{i j}^{*} V_{i} V_{j}^{*} \quad(i, j) \in E \cup F 
\end{align}
\begin{align}
\label{10h}
-\theta_{i j}^{\Delta} \leqslant \angle\left(V_{i} V_{j}^{*}\right) \leqslant \theta_{i j}^{\Delta} \quad\forall(i, j) \in E
\end{align}
\begin{align}
\label{10i}
\tan \left(-\theta_{i j}^{\Delta}\right) \Re\left(V_{i} V_{j}^{*}\right) \leqslant \Im\left(V_{i} V_{j}^{*}\right) \leqslant \tan \left(\theta_{i j}^{\Delta}\right) \Re\left(V_{i} V_{j}^{*}\right)
\end{align}
\end{subequations}

From above equations, it can be realized that the non-convex nature of the OPF problem arises from product of the voltage variables $V_{i} V_{j}^{*}$. Assume, a new lifted variable $W_{ij}$ in place of $V_{i} V_{j}^{*}$ for applying convex relaxation methods.

\begin{align}
\label{eq:W_{i j}=V_{i} V_{j}^{*}}
W_{i j}=V_{i} V_{j}^{*}
\end{align}

Equations~\eqref{10g}, \eqref{10h}, and~\eqref{10i} can be written as,

\begin{subequations}
\begin{align}
&S_{i j}=Y_{i j}^{*} W_{i j}-Y_{i j}^{*} W_{i j} \quad(i, j) \in E \cup F\\ 
&-\theta_{i j}^{\Delta} \leqslant \angle\left(W_{i j}\right) \leqslant \theta_{i j}^{\Delta} \quad\forall(i, j) \in E\\
&\tan \left(-\theta_{i j}^{\Delta}\right) \Re\left(W_{i j}\right) \leqslant \Im\left(W_{i j}\right) \leqslant \tan \left(\theta_{i j}^{\Delta}\right) \Re\left(W_{i j}\right)
\end{align}
\label{eq:relax1}
\end{subequations}

Next we explain how defining lifted variables can be leveraged to convexify the OPF problem.  

\section{Convex Relaxations in Optimal Power flow}
\label{Convex Relaxations in Optimal Power flow}
Traditional OPF solving methods may find global optima of the solution but they might stuck in local optima~\cite{molzahn2013sufficient}. Conversely, convex relaxation techniques can obtain bounds on the optimal objective values, certify infeasibility, and in some cases, achieve globally optimal solutions. There has been many relaxation techniques applied to the OPF problems i.e Second order cone relaxations, Quadratic Convex relaxations and Semi-definite relaxations~\cite{low2014convex},~\cite{molzahn2013implementation}. Each relaxation method follows separate methodology to solve the non-convex OPF problem. To be specific, they approach differently to convexify the source of non-convexification $V_{i} V_{j}^{*}$. In this study we will focus only on quadratic convex (QC) relaxation. 

\textbf{Quadratic Convex Relaxations}: The quadratic convex (QC) relaxation is a approach that encloses the trigonometric and product terms in the polar representation of power flow equations within convex envelopes~\cite{Narimani3,Narimani4,narimani2020tightening}. It represents the voltage variables in polar coordinates and expands equation~\eqref{eq:W_{i j}=V_{i} V_{j}^{*}} in following way,
\begin{subequations}
\label{eq:qc convex}
\begin{align}
\label{eq:objective}
 W_{i i} = v_i^2\quad \forall i \in N 
\end{align}
\begin{align}
\Re{(W_{ij})} = {v_i}{v_j}\cos({\theta{i}-\theta{j})}\quad \forall(i, j) \in E
\end{align}
\begin{align}
\Im{(W_{ij})} = {v_i}{v_j}\sin({\theta{i}-\theta{j})}\quad \forall(i, j) \in E
\end{align}
\end{subequations}
QC relaxation relaxes these constraints by drawing tight envelopes around nonconvex terms. Such as, convex envelopes for square terms can be defined as~\cite{mccormick1976computability},
\begin{align}
\label{eq:convex square}
\langle x^{2}\rangle^{T} \equiv \begin{cases}
\widecheck{x} \geqslant x^{2} \\
\widecheck{x} \leqslant\left(x^u+x^l\right) x-x^u x^l.\\
\end{cases}
\end{align}
Here, ${\widecheck{x}}$, $x^u$ and $x^l$ corresponds to convex envelopes of $x$, upper and lower bound of $x$, respectively.
Additionally, convex envelopes for bilinear terms can be defined as,

\begin{align}
\label{eq:convex product}
\langle x y\rangle^{M} \equiv\left\{\begin{array}{l}
\widecheck{x y} \geqslant x^{l} y+y^{l} x-x^{l} y^{l} \\
\widecheck{x y} \geqslant x^{u} y+y^{u} x-x^{u} y^{u} \\
\widecheck{x y} \leqslant x^{l} y+y^{u} x-x^{l} y^{u} \\
\widecheck{x y} \leqslant x^{u} y+y^{l} x-x^{u} y^{l}
\end{array}\right.
\end{align}

Convex envelopes for sine and cosine function for $x \in [0,\pi/2]$ can be given as~\cite{hijazi2017convex},

\begin{align}
\begin{array}{c}
\label{eq:trig}
\langle\sin (x)\rangle^{S} \equiv\left\{\begin{array}{l}
\widecheck{S} \leqslant \cos \left(\frac{x^{u}}{2}\right)\left(x-\frac{x^{u}}{2}\right)+\sin \left(\frac{x^{u}}{2}\right) \\
\widecheck{S} \geqslant \cos \left(\frac{x^{u}}{2}\right)\left(x+\frac{x^{u}}{2}\right)-\sin \left(\frac{x^{u}}{2}\right)
\end{array}\right. \\
\langle\cos (x)\rangle^{C} \equiv\left\{\begin{array}{l}
\widecheck{C} \leqslant 1-\frac{1-\cos \left(x^{u}\right)}{\left(x^{u}\right)^{2}} x^{2} \\
\widecheck{C} \geqslant \cos \left(x^{u}\right)
\end{array}\right.
\end{array}
\end{align}

Now using equations \eqref{eq:qc convex}--\eqref{eq:trig} convex relaxations of product terms in power flow equation can be obtained as follows.

\begin{subequations}
\label{eq:non-convex}
\begin{align}
&W_{i i}=\left\langle v_{i}^{2}\right\rangle^{T} \quad i \in N \\
&\Re\left(W_{i j}\right)=\left\langle\left\langle v_{i} v_{j}\right\rangle^{M}\left\langle\cos \left(\theta_{i}-\theta_{j}\right)\right\rangle^{C}\right\rangle^{M} \quad \forall(i, j) \in E \\
&\Im\left(W_{i j}\right)=\left\langle\left\langle v_{i} v_{j}\right\rangle^{M}\left\langle\sin \left(\theta_{i}-\theta_{j}\right)\right\rangle^{S}\right\rangle^{M} \quad \forall(i, j) \in E
\end{align}
\end{subequations}

Incorporating equations~\eqref{eq:relax1}--\eqref{eq:non-convex} into equation \eqref{OPF Formulation overall} would convexify the nonconvex terms and result in the QC relaxation of the OPF problem. The solution of the QC relaxation is a lower bound for the original nonconvex OPF problem. The tighter the relaxation, the better a lower bound can be computed for the OPF problem. Next, we leverage the graphlets analysis to understand the correlation between power system topology and the optimality of the QC relaxation for different test cases.



\section{Network Graphlet}
\label{Network Motifs}
\begin{figure}[H]
    \centering
    \captionsetup{justification=centering}
    \subfloat[\centering Graphlet type 1 ]{{\includegraphics[width=2cm]{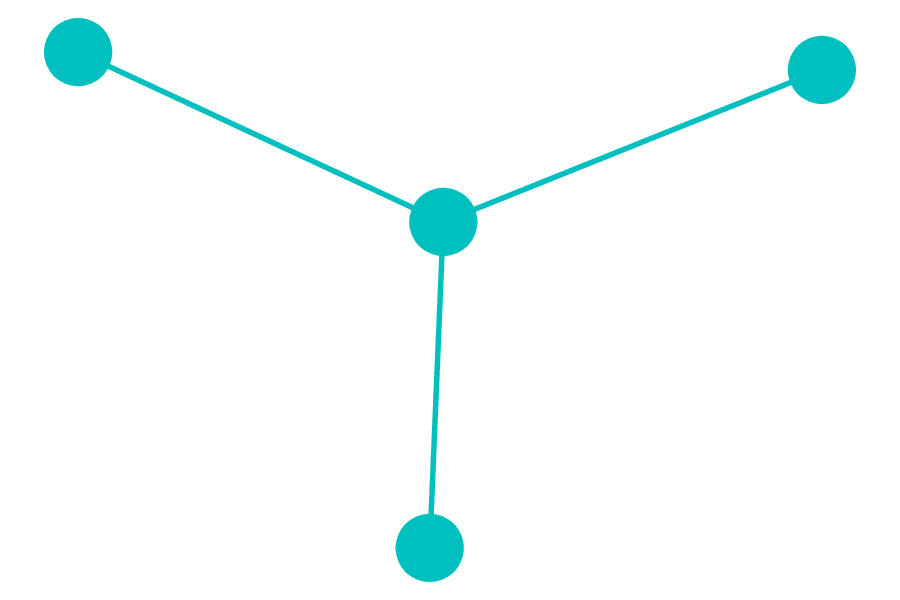} }}
    \qquad
    \subfloat[\centering Graphlet type 2]{{\includegraphics[width=2cm]{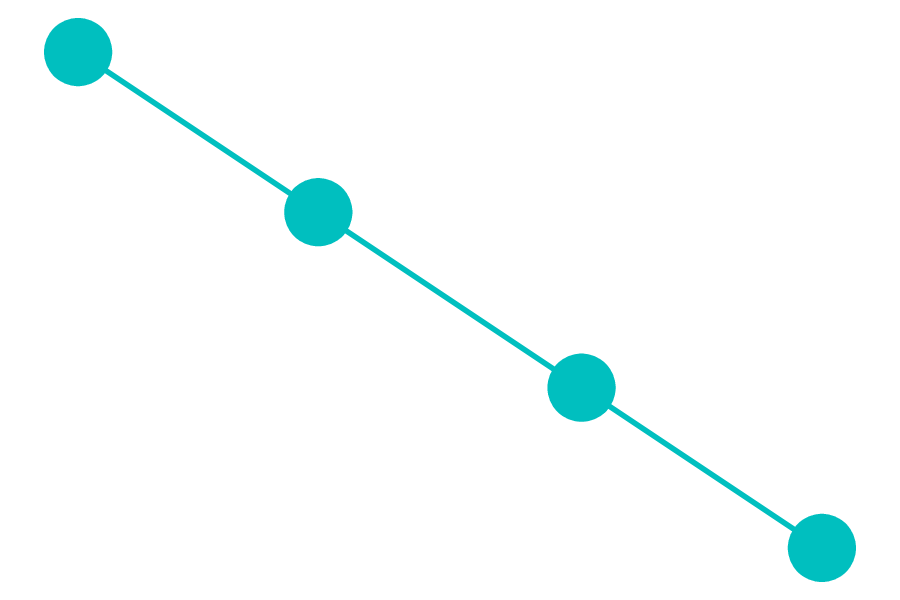} }}%
    \subfloat[\centering Graphlet type 3]{{\includegraphics[width=2cm]{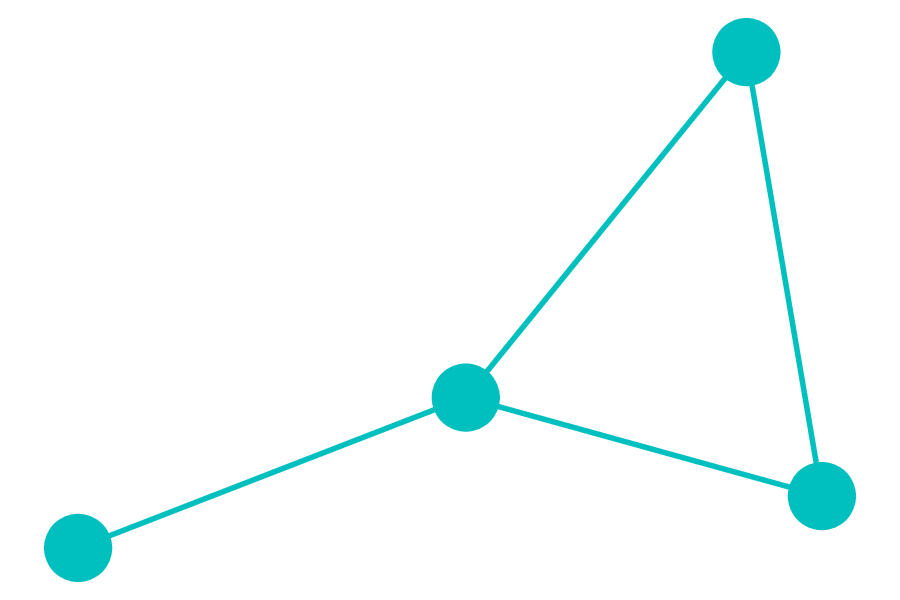} }}
    \qquad
    \subfloat[\centering Graphlet type 4]{{\includegraphics[width=2cm]{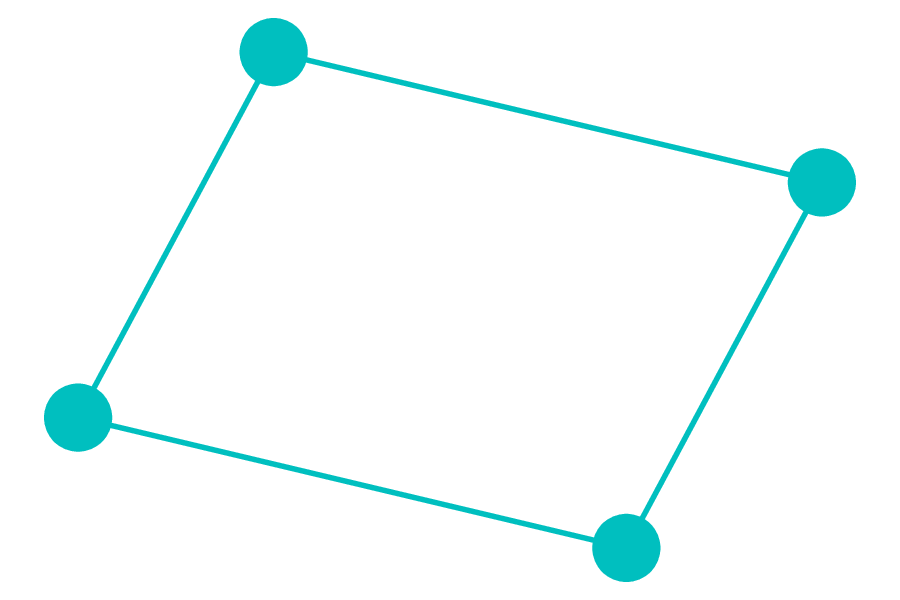} }}
    \qquad 
    \subfloat[\centering Graphlet type 5]{{\includegraphics[width=2cm]{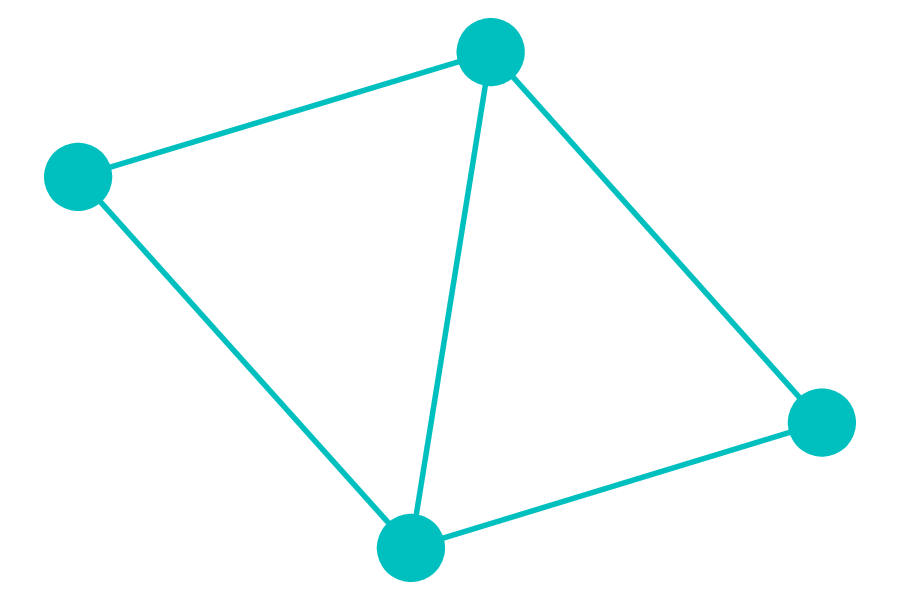} }}
    \qquad
    \subfloat[\centering Graphlet type 6]{{\includegraphics[width=2cm]{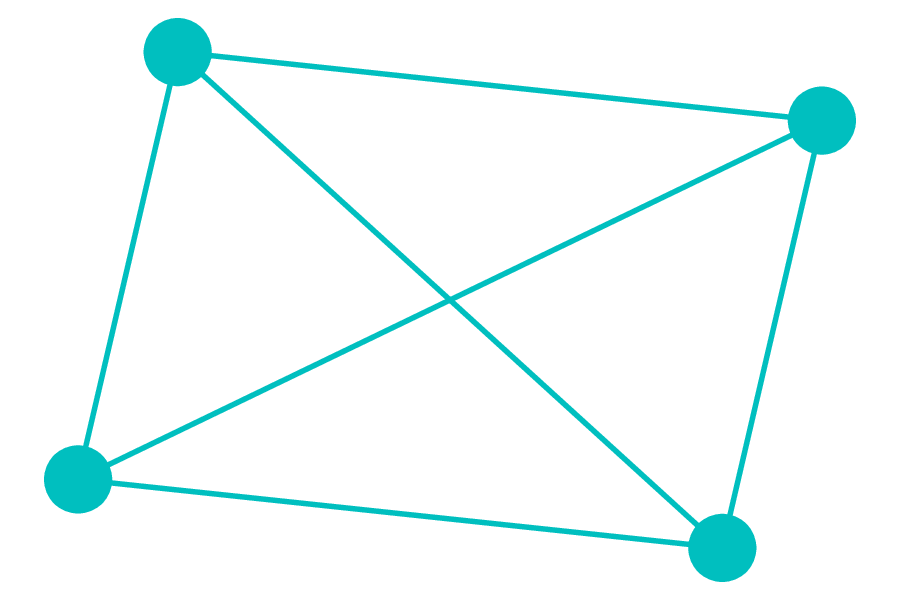} }}
    \qquad
    \caption{Different graphlet types}
    \label{fig:different motifs}
\end{figure}

To observe the local structure of any power grid we consider the grid as a undirected graph $G(V,E)$ where $V$ as nodes is the buses in the grid and $E$ as edges stands for Lines in the grid. The order and the size of graph $G$ can be defined as total number of nodes and total number of edges of the graph. A graph $G'(V',E')$ is a subgraph of graph ${G}$ if ${G'} \subseteq {G}\;  \text{such that}\; {V'} \subseteq {V}\;{and}\; {E'} \subseteq {E}$. That subgraph $G'$ can be called as induced subgraph of $G$ if $E'$ contains all edges $e_{uv}\in E$ such that $u,v\in V'$. Graphs $G'$ and $G''$ can be called isomorphic if there exists a bijection $h:V'\rightarrow V''$ such that any two adjacent nodes $u,v \in V' $ of $G'$ are also adjacent in $G''$ after the mapping occurs. If $G_k = (V_k,E_k)$ is a $k$ node subgraph of G and there exists an isomorphism between $G_k$ and $G'$ where $G' \in G$, then there is an occurrence of $G_k$ in $G$. A motif can be defined as a multi node subgraph pattern that occurs too frequently in a graph comparing to some random graphs. These recurrent patterns are considered as building blocks of networks, and different combinations of a small number of motifs can generate enormously diverse forms. Albeit the notion of motif originated in biological networks, different motifs could be found in a variety of different types of complex networks. Traditional network attributes, which characterize nodes, connections, or the entire network, do not really afford for the profiling of local structural characteristics. As it is well known, motif patterns differ substantially among networks~\cite{xia2019survey}. In this paper, we use 4 node connected undirected subgraphs which are called graphlets. Therefore, 6 types of graphlets that can be found in the power network are used to asses the local structure. These graphlets are shown in figure~\ref{fig:different motifs}. In this paper, we first explore the existence of different graphlet types among different buses of the network and then classify those subgraph patterns.

\section{Detection of Graphlets} 
\label{Statistical Analysis}
In this section, how graphlet can be identified in different networks is discussed.

\begin{figure}[H]
\centering
\captionsetup{justification=centering,margin=6cm}
{\includegraphics[width=4.4cm]{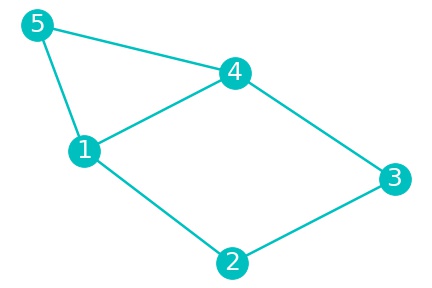}}
\caption{A network having 5 nodes and 6 edges. Total of 5 graphlets can be found in this network}
\label{fig:type 3 motif}
\end{figure}

This paper uses an algorithm that follows the footprints in~\cite{wernicke2006fanmod}\cite{ribeiro2009strategies}.
Given a graph G = (V, E), the following algorithm enumerates all of its 4-nodes subgraphs. There is two sets, one is subgraph set, and another is extension set. Extension set consists of neighboring nodes that are not in subgraph set. When summation of subgraph set and extension set element equals to our desired subgraph level which is four, this algorithm returns a subgraph.
\begin{itemize}
  \item In the first stage, the subgraph level is one which is basically a node. In figure \ref{fig:type 3 motif}, for this stage, we consider node 1 as a subgraph set. Therefore, extension set consists of \{2, 4, 5\}. For node 2 as a subgraph set, extension set consists of node {3}. However, Node 1 is a neighboring node of node 2 but is not considered in extension set as it has already been considered as a subgraph set. Following this formula, other combinations of subgraph and extension set are [\{2\},\{3\}], [\{3\},\{4\}], [\{4\},\{5\}]. This stage returns only one 4 node subgraph, \{{1,2,4,5}\}.
  \item In the second stage, the subgraph level is two which is an edge. In this stage, for node \{1, 2\} as a subgraph set, extension set consists of neighboring nodes of edge 1-2 which is \{3, 4, and 5\}. For edge \{1,4\} as a subgraph set, extension set consists of \{3,5\}. Here node 2 is not in extension set as edge {1,2} already have been considered as a subgraph set. Only other possible combination are [\{2,3\},\{4,5\}]. This stage returns two 4 node subgraph. Those are, 
  \{1,4,3,5\} and \{2,3,4,5\}.
  \item Following this formula for third and fourth stage, this algorithm gives two more subgraph combination. Those are \{1,2,3,4\} and \{1,2,3,5\}.
  \item After that, this algorithm classifies those subgraphs or graphlets. For example, \{2, 3, 4, 5\} is a type 2 graphlet as previously mentioned. Following the classification, occurrence of each graphlet is counted. Such as, type 2 graphlet has occurred 4 times out of total 5 graphlets in this graph. Furthermore, for all other graphlets its graphlet occurrence time divided by total graphlet count gives us a percentage. For different test cases in Matpower~\cite{zimmerman1997matpower} and PGLib-OPF v18.08 benchmark library~\cite{babaeinejadsarookolaee2019power}, this percentages for all graphlet types are given in Table \ref{table:motif percnt all cases}.
\end{itemize}

\begin{table*}[!ht]
    \centering
    \caption{Different Networks with different graphlet percentages and their corresponding QC relaxation optimality gap}
    \label{table:motif percentage all}
    \begin{tabular}{|l|c|c|c|c|c|c|c|c|}\hline
        \thead{Networks} & \thead{QC Gap} & \thead{Type 1} & \thead{Type 2} & \thead{Type 3} & \thead{Type 4} & \thead{Type 5} & \thead{Type 6} \tabularnewline \hline
        case18 & 1.17E-05 & \textbf{15.38461538} & 84.61538462 & 0 & 0 & \textbf{0} & \textbf{0} \tabularnewline \hline
        case9 & 0.000379748 & \textbf{20} & 80 & 0 & 0 & \textbf{0} & \textbf{0} \tabularnewline \hline
        pglib\_opf\_case\_ACTIVSg200 & 0.003350142 & 39.50471698 & 54.54009434 & 5.247641509 & 0.471698113 & 0.235849057 & \textbf{0} \tabularnewline \hline
        pglib\_opf\_case24\_ieee\_rts & 0.01203402 & \textbf{21.51162791} & 72.6744186 & 2.325581395 & 3.488372093 & \textbf{0} & \textbf{0} \tabularnewline \hline
        pglib\_opf\_case39\_epri & 0.02275154 & \textbf{21.29032258} & 75.48387097 & 1.935483871 & 1.290322581 & \textbf{0} & \textbf{0} \tabularnewline \hline
        pglib\_opf\_case73\_ieee\_rts & 0.030513151 & \textbf{21.25603865} & 73.2689211 & 2.576489533 & 2.898550725 & \textbf{0} & \textbf{0} \tabularnewline \hline
        pglib\_opf\_case30\_as & 0.055670456 & \textbf{22.07207207} & 60.36036036 & 16.66666667 & 0.900900901 & \textbf{0} & \textbf{0 }\tabularnewline \hline
        pglib\_opf\_case60\_c & 0.058817246 & \textbf{24.46555819} & 68.64608076 & 6.413301663 & 0.237529691 & 0.237529691 & \textbf{0} \tabularnewline \hline
        pglib\_opf\_case2868\_rte & 0.094801158 & 37.85541316 & 58.46053383 & 3.202459663 & 0.440780344 & 0.040812995 & \textbf{0} \tabularnewline \hline
        pglib\_opf\_case14\_ieee & 0.10907652 & 10.14492754 & 62.31884058 & 24.63768116 & 0 & 2.898550725 & \textbf{0} \tabularnewline \hline
        pglib\_opf\_case2848\_rte & 0.118583008 & 36.59196075 & 59.56681719 & 3.275352623 & 0.521268886 & 0.044600546 & \textbf{0} \tabularnewline \hline
        pglib\_opf\_case1951\_rte & 0.126574378 & 38.41369779 & 57.3786855 & 3.62791769 & 0.525952088 & 0.053746929 & \textbf{0} \tabularnewline \hline
        pglib\_opf\_case179\_goc & 0.151274263 & 29.96515679 & 59.93031359 & 8.797909408 & 0.609756098 & 0.609756098 & 0.087108014 \tabularnewline \hline
        pglib\_opf\_case57\_ieee & 0.158611501 & 15.15151515 & 69.6969697 & 14.04958678 & 0.550964187 & 0.550964187 & 0 \tabularnewline \hline
        pglib\_opf\_case89pegase & 0.165595479 & 14.37705201 & 43.13115604 & 30.43027475 & 0.380162433 & 7.033005011 & 4.648349749 \tabularnewline \hline
        pglib\_opf\_case500\_goc & 0.240470483 & 32.58928571 & 59.56909938 & 7.104037267 & 0.36878882 & 0.349378882 & 0.019409938 \tabularnewline \hline
        pglib\_opf\_case2737sop\_k & 0.256111524 & 27.69497108 & 69.78001098 & 2.068994638 & 0.434911118 & 0.02111219 & 0 \tabularnewline \hline
        pglib\_opf\_case2000\_goc & 0.304646737 & 23.75966093 & 68.35203191 & 6.881077038 & 0.513587634 & 0.463724757 & 0.029917726 \tabularnewline \hline
        pglib\_opf\_case2746wp\_k & 0.313345858 & 27.58125185 & 69.87870335 & 2.096276573 & 0.422636406 & 0.02113182 & 0 \tabularnewline \hline
        pglib\_opf\_case3970\_goc & 0.337410369 & 20.84741502 & 76.12211267 & 2.834957964 & 0.19272128 & 0.002793062 & 0 \tabularnewline \hline
        pglib\_opf\_case3375wp\_k & 0.530375619 & 27.47302905 & 68.91286307 & 2.933609959 & 0.257261411 & 0.315352697 & 0.107883817 \tabularnewline \hline
        pglib\_opf\_case9591\_goc & 0.607902624 & 21.44041766 & 74.65971875 & 3.599572134 & 0.280664567 & 0.017664203 & 0.001962689 \tabularnewline \hline
        pglib\_opf\_case118\_ieee & 0.785063339 & \textbf{26.102169 }& 62.841148 & 9.237229 & 1.39958 & 0.349895 & 0.069979 \tabularnewline \hline
        pglib\_opf\_case2853\_sdet & 0.863872504 & \textbf{40.43575845} & 44.6747692 & 11.9717117 & 0.256203947 & 0.741085797 & 1.920470907 \tabularnewline \hline
        pglib\_opf\_case793\_goc & 1.316740725 & \textbf{28.8169364}9 & 69.4718876 & 1.320049813 & 0.149439601 & 0.231416564 & 0.027687821 \tabularnewline \hline
        pglib\_opf\_case6468\_rte & 1.744539041 & \textbf{37.27092714} & 58.57706812 & 3.613454917 & 0.469181581 & 0.06306204 & 0.006306204 \tabularnewline \hline
        pglib\_opf\_case4917\_goc & 2.491041794 & \textbf{38.3252531} & 57.84614288 & 3.156637673 & 0.50798517 & 0.144374733 & 0.019606445 \tabularnewline \hline
        pglib\_opf\_case240\_pserc & 2.727162682 & 23.77777778 & 63.14814815 & 11.11111111 & 0.925925926 & 1 & 0.037037037 \tabularnewline \hline
        pglib\_opf\_case3022\_goc & 2.752412865 & \textbf{41.4685445} & 53.63541972 & 4.068049306 & 0.48454514 & 0.279545273 & 0.063896062 \tabularnewline \hline
        pglib\_opf\_case162\_ieee\_dtc & 5.838058785 & 17.08751297 & 65.0639917 & 14.07817364 & 0.864752681 & 2.525077828 & 0.38049118 \tabularnewline \hline
        pglib\_opf\_case6515\_rte & 6.389830911 & \textbf{37.67528825} & 58.22000623 & 3.573698972 & 0.462449361 & 0.062324712 & 0.006232471 \tabularnewline \hline
        pglib\_opf\_case6495\_rte & 15.08381955 & \textbf{37.42583467} & 58.43473451 & 3.603932019 & 0.466361625 & 0.062851971 & 0.006285197 \tabularnewline \hline
    \end{tabular}
    \label{table:motif percnt all cases}
\end{table*}

\section{Results}
\label{Results}

We implement the QC relaxation and local structure of power systems' graph on various test cases from PGLib-OPF v18.08 benchmark library~\cite{babaeinejadsarookolaee2019power} and Matpower to find the correlation between local structure of power systems' graph and OPF optimality gap. 
Table~\ref{table:motif percnt all cases} tabulates optimality gap from QC relaxation and their corresponding different graphlet type percentages for different test cases. This table is sorted with  ascending "QC gap" values. Clearly, some patterns are visible in this Table. The most noticeable trend is that, of the 33 cases investigated, 12 cases with the low optimality gap have one characteristic in common. That is they share the common trait of having no type 6 graphlets. From figure \ref{fig:different motifs}, it is apparent that type 6 graphlet is equivalent to 4-nodes complete subgraph. However, existence of complete 4-node subgraph is not very frequent in any kind of networks. That is why in other cases such low percentage of type 6 graphlets are observed. For instance, "pglib\_opf\_case162\_ieee\_dtc" has total 2891 graphlets. Yet, it has only 11 type 6 graphlets. Particularly, this trend shows us that existence of type 6 graphlets even in such low percentage can affect the QC relaxation optimality gap. Additionally, there is another pattern that highlights less percentage of type 1 graphlets in low optimality gap networks in comparison with networks with high optimality gaps. Particularly, networks with low optimality gap have 15-25\%  type 1 graphlets whereas most networks having high optimality gap have 25-40\% type 1 graphlets. Furthermore, not as significant pattern as type 6, type 5 graphlets are also absent in some of the lowest optimality gap networks. In conclusion, from Table \ref{table:motif percnt all cases}, it can be observed that existence of type 6, type 1 and type 5 graphlets can significantly impact OPF optimality gap in different test cases.

The graphlet structure can be leveraged to introduce constraints or relationships between variables in the convex relaxation problems. For instance, in a graph with specific graphlets, the variables within the graphlets may have strong dependencies or correlations. By incorporating these graph constraints into the relaxation formulation, the relaxation can better capture the underlying structure and improve its tightness. Moreover, the graph structure can provide insights into the connectivity of variables. If certain subsets of variables are strongly connected in the graph, it suggests that they should have similar values or follow certain patterns. By enforcing such connectivity constraints in the relaxation, the relaxation can better capture the relationships between variables and improve its accuracy. Overall, by incorporating and leveraging the graph structure in the formulation of the relaxation, it is possible to exploit the inherent properties of the problem and improve the tightness of the relaxation. The proposed approach results in more accurate solutions and better bounds for convex relaxation.

\section{Conclusion}
\label{Conclusion}

This paper investigates the correlation between the optimality gap of the OPF convex relaxation and the local structure of power system graphs addressed by the distribution of graphlet types. In this context, the paper calculates the percentage of different 4-node graphlets in various networks and examines the relationship between graphlet counts and optimality gaps for the QC relaxation of the OPF problem. Consequently, the results clearly indicate that networks with high QC optimality gaps share a common characteristic of having type 6 graphlets or complete 4-node graphlets. Additionally, it is noticeable that networks with high QC optimality gaps tend to have a high percentage of type 1 or 4-node star-shaped graphlets. Conversely, most networks with very low optimality gaps do not have type 5 graphlets. In conclusion, this study suggests that the existence of certain graphlets can cause a high optimality gap in the QC relaxation. In the future, this research can be expanded to associate electrical parameters, such as branch admittance, with graphlet-level analysis. This research is currently studying the identification of nodes in the power system where enforcing redundant constraints can tighten the QC relaxation of the OPF problem for those test cases.

Regenerate response 
\bibliographystyle{IEEEtran}
\bibliography{ref}
\end{document}